\documentclass{amsart}
\usepackage{amssymb}
\usepackage{amsmath}
\usepackage{amscd}
\usepackage{bbm}
\usepackage{accents}
\usepackage[all]{xy}

\usepackage{ucs}
\usepackage[utf8x]{inputenc}
\usepackage[T1,T2A]{fontenc}
\usepackage[english]{babel}

\newcommand\dd{\partial}

\newcommand\Rs{{\mathbb R}}

\newcommand\ddiv{\mathop {\fam 0 div}\nolimits} 
\newcommand\curl{\mathop {\fam 0 curl}\nolimits}

\begin{document}

\title{A steady smooth Euler flow with support in the  vicinity of  a helix}


\author{A.V.Gavrilov } 

\email{gavrilov19@gmail.com}

\maketitle 

\begin{abstract} 
In this article we construct a smooth Euler flow supported in a  neighborhood of a helix. It may be considered  a generalization of a similar solution found  by the author  
for a circle. 
\end{abstract}

\section{Introduction}

In this article we construct a smooth Euler flow in $\Rs^3$ supported in a  neighborhood of a helix. While it 
may be considered  a generalization of the  solution of the Euler equation  found   by the author  in \cite{Gav}, it is not of very much  interest by itself. 
(The whole point of \cite{Gav}  was to find a smooth Euler flow with compact support. Of course, a helix  is completely  useless for this purpose.) 
Maybe  this solution  could be a stepping-stone  to  more interesting  generalizations.   

It is convenient to interpret  this new flow as a  modification of the old one, and  for this reason the notation  we use here is  similar to \cite{Gav}.
We consider a helix ${\mathcal C}\subset\Rs^3$ described in the standard cylindrical coordinates $(\rho, \varphi, z)$  by the equations 
$\rho=1,\,z=k\varphi$ where $k>0$ is the slope. This curve has  a group of  isometries  with the generator   
$\xi=\dd_\varphi+k\dd_z$, taking advantage of this  we construct a flow which is invariant under this isometries.  
This flow retains simple topology described by  the Arnold's theorem \cite[Ch. II, Theorem 1.2]{AK}  except invariant tori become invariant cylinders.

\section{The Euler equation in cylindrical coordinates}

\subsection{Preliminaries}

A vector field in cylindrical coordinates  is usually  written 
using the local  basis   $e_\rho=\dd_\rho,\,e_z=\dd_z,\,e_\varphi=\frac{1}{\rho}\dd_\varphi$.
In this coordinates  the incompressibility condition $\ddiv u=0$  for  $u=u_\rho e_\rho+u_z e_z+u_\varphi e_\varphi $ becomes 
$$\ddiv  u=\dd_\rho u_\rho+\frac{1}{\rho}u_\rho+\frac{1}{\rho}\dd_\varphi u_\varphi+\dd_z u_z=0,\eqno{(1)}$$
and the Euler equation itself $(u\cdot\nabla) u=-\nabla p$ is \cite[\S 15]{LL} 
\begin{equation*}
\left\{
\begin{array}{rl}
u_\rho\frac{\dd}{\dd\rho}u_\rho+\frac{1}{\rho}u_\varphi\frac{\dd}{\dd\varphi}u_\rho+u_z \frac{\dd}{\dd z}u_\rho-\frac{1}{\rho}u_\varphi^2=-\frac{\dd}{\dd\rho}p, \,\,\,\,\,\,\,\,\,\, (2a)\\
u_\rho\frac{\dd}{\dd\rho}u_\varphi+\frac{1}{\rho}u_\varphi\frac{\dd}{\dd\varphi}u_\varphi+u_z \frac{\dd}{\dd z}u_\varphi+\frac{1}{\rho}u_\rho u_\varphi=-\frac{1}{\rho}\frac{\dd}{\dd\varphi}p,\,\,\,\,\,\,\,\,\,\, (2b)\\
u_\rho\frac{\dd}{\dd\rho}u_z+\frac{1}{\rho}u_\varphi\frac{\dd}{\dd\varphi}u_z+u_z \frac{\dd}{\dd z}u_z=-\frac{\dd}{\dd z}p. \,\,\,\,\,\,\,\,\,\, (2c)
\end{array} \right. 
\end{equation*}

Following \cite{Gav}, we assume that 
$$|u|^2=u_\rho^2+u_z^2+u_\varphi^2=3p. \eqno{(3)}$$
(Note that with this additional equation the system (1-3) becomes  overdetermined.) In this case the Bernoulli law
$$u\cdot\nabla\left(\frac{1}{2}|u|^2+p\right)=0$$
implies  $u\cdot\nabla p=0$.  In general, we will call an Euler flow satisfying the latter condition \emph{localizable}
(for reasons explained in \S 4.2).

\subsection{The flow}  

We are looking  for a solution of the form 
$$u=\frac{1}{x}\frac{\dd t}{\dd y} e_\rho+\frac{1}{x^2+k^2}\left(kh-x\frac{\dd t}{\dd x} \right)e_z+\frac{1}{x^2+k^2}\left(xh+k\frac{\dd t}{\dd x} \right)e_\varphi,\eqno{(4)}$$
where $t$ and $h$ are  functions of the variables
$$x=\rho,\,y=z-k\varphi.$$
This field is obviously invariant under the isometries, and it is not difficult to check that $\ddiv u=0$. 
Also, it follows directly from (4) that $u\cdot\nabla t=0$, hence  to satisfy one of the necessary conditions  $u\cdot\nabla p=0$ it is sufficient to assume\footnote{In fact, this assumption is  more or less  unavoidable.}
that $p=p(t)$. There is some leeway in  chosing this function (due to the modification discussed in Sec. 4.2), a convenient  choice  is 
$$p=\frac{t}{1+k^2}.\eqno{(5)}$$

\subsection{The equations}

Combining (2b) and (2c) we have 
$$\rho(u_\rho\frac{\dd}{\dd\rho}u_\varphi+\frac{1}{\rho}u_\varphi\frac{\dd}{\dd\varphi}u_\varphi+u_z \frac{\dd}{\dd z}u_\varphi+\frac{1}{\rho}u_\rho u_\varphi)+k(u_\rho\frac{\dd}{\dd\rho}u_z+\frac{1}{\rho}u_\varphi\frac{\dd}{\dd\varphi}u_z+u_z \frac{\dd}{\dd z}u_z)=-\xi(p)=0,$$
or simply 
$$u\cdot\nabla h=0.$$
This equation would follow if we assume\footnote{This is also not very much of an assumption because 
it is not difficult to see  that $dt\wedge dh=0$.} that $h=h(t)$ is also a function of $t$.   
$$$$
 What is left is the remaining pair of the Euler equations  together with (3), which takes the form 
$$\frac{1}{x^2}\left(\frac{\dd t}{\dd y}\right)^2+\frac{1}{x^2+k^2}\left[\left(\frac{\dd t}{\dd x}\right)^2+h^2\right]=\frac{3t}{1+k^2}. \eqno{(6)}$$
Changing the variables in (2a), (2c) from $(\rho, z)$ to $(x,y)$ we then have 
$$(t_yt_{xy}-t_xt_{yy})-\frac{t_y^2}{x}-x\left(\frac{xh+kt_x}{x^2+k^2}\right)^2+\frac{x^2t_x}{1+k^2}=0,  \eqno{(7a)}$$
$$\left(t_yt_{xx}-t_xt_{xy}\right)-\frac{t_xt_y(x^2-k^2)}{x(x^2+k^2)}+\frac{2kht_y}{x^2+k^2}-\frac{(x^2+k^2)t_y}{1+k^2}=0.\eqno{(7b)}$$

\section{The solution}

\subsection{The change of variables}

Following the same tactic  as in \cite{Gav}  we  assume that   
$$\frac{\dd}{\dd x}t=F,\,\left(\frac{\dd}{\dd y}t\right)^2=G$$ 
where $F,G$ are  functions of $(x,t)$. This  assumption obviously  implies that this two functions  satisfy the following partial differential  equation 
$$\frac{\dd G}{\dd x}+F\frac{\dd G}{\dd t}=2G\frac{\dd F}{\dd t}.\eqno{(8)}$$ 
What we want is to rewrite  the equations (6-7) in terms of  $F$ and $G$.   To begin with, (6) simply  turns into an algebraic relation
$$G=x^2\left(\frac{3t}{1+k^2}-\frac{F^2+h^2}{x^2+k^2}\right). \eqno{(9)}$$
Taking  into account that $t_x=F,\,t_{xx}=F_x+FF_t,\,\,t_{xy}=t_yF_t,$
we may write  (7b)  as  
$$F_x-\frac{x^2-k^2}{x(x^2+k^2)}F=-\frac{2kh}{x^2+k^2}+\frac{x^2+k^2}{1+k^2}.\eqno{(10)}$$
This linear differential  equation has a solution
$$F=\frac{kh}{x}+\frac{(x^2+k^2)(x^2-c)}{2x(1+k^2)}\eqno{(11)}$$
where  $c=c(t)$. Finally, using $t_y^2=G,\,\,t_{yy}=\frac{1}{2}G_t$ and (8), the last equation (7a) may be rewritten as 
$$\frac{\dd}{\dd x}\left(\frac{G}{x^2}\right)=-\frac{2F}{1+k^2}+\frac{2}{x} \left(\frac{xh+kF}{x^2+k^2} \right)^2,$$
which is actually  a consequence of (9) and (10). 

\subsection{The ODE}

Under the assumptions we have made all  the original equations (1-3) are satisfied. However, there is also the new one  (8) 
which  is not done yet. This equation contains two unknown functions of  $t$, namely $h$ and $c$.   After  substituting (9) and  (11)  into (8) 
and obvious  algebraic transformations, it is possible to get rid of  the variable $x$  and reduce this  PDE to  two (rather cumbersome) ordinary differential  equations, 
$$\frac{dh}{dt}=\frac{(k^2+c)S+6t(k^2+1)(kh+6t)}{2(1+k^2)(hS+18kt^2)},\,\,\,\frac{dc}{dt}=\frac{kS+6th(k^2+1)}{hS+18kt^2},\eqno{(12)}$$
where 
$$S=h^2(1+k^2)-3t(c+k^2).$$
We are interested in a solution of this system with initial condition 
$$h(0)=0,\,\,c(0)=1.$$
Note that the denominator at this point becomes  zero, so this is a singular Cauchy problem. It has no analytic  solutions, but one can show\footnote{Apparently, it is not possible to reduce this system to a Briot-Bouquet equation the way it is done in \cite{Gav}. We   have to  prove this fact using  the  series directly,  which  is   straightforward but somewhat bothersome.} 
that it has a solution in the form of a  Puiseux series, analytic as a function of  $s=\sqrt{t}$ 
$$h=s+\frac{8k}{3(1+k^2)}s^2+\frac{82k^2-189}{36(1+k^2)^2}s^3+\frac{k(154k^2-369)}{27(1+k^2)^3}s^4+O(s^5),$$ 
$$c=1+2ks+\frac{10k^2-9}{3(1+k^2)}s^2+\frac{k(106k^2-249)}{18(1+k^2)^2}s^3+\frac{3688k^4-8658k^2+243}{216(1+k^2)^3}s^4+O(s^5).\eqno{(13)}$$

\section{Completing the construction}

\subsection{The variable   $y$}

Now we have to change the variables from $(x,t)$ back to $(x,y)$. This part is slightly more complicated then in \cite{Gav}
because in this case we have \emph{two} completely different solutions instead of just one. It is easy to see why this happens 
if we take a closer look at the geometry. The streamlines of the original flow $u$ in \cite{Gav} have the form of slightly deformed 
helices winding around the circle ${\mathcal C}$. We still have the same picture when ${\mathcal C}$ itself becomes a helix, 
except in this case it does matter if the helicity of ``small'' helices is the same as the ``big'' one or the opposite.  
The first choice  corresponds to $s>0$, and the second one to $s<0$. 

The function $G$ given by (9) is obviously real analytic as a function of $x$ and $s=\pm\sqrt{t}$ at the point $(x,s)=(1,0)$.
However, as it is supposed to be equal to the  square of $\frac{\dd t}{\dd y}$, the region $G<0$  is forbidden.
A direct computation shows that the condition $G\ge 0$ is equivalent to 
$$x^6+(k^2-2c)x^4+[4(1+k^2)(kh-3t)+c^2-2k^2c]x^2+4(1+k^2)(hk^2-kc+h)h+k^2c^2\le 0,\eqno{(14)}$$
which means (using  $X=x-1$)   
$$4(1+k^2)(X^2-2s^2)+4(3+k^2)X^3-16kX^2s-4(2k^2+3)Xs^2+\frac{64}{3}ks^3\le O(X^4+s^4).$$
This domain consists of two  parts corresponding to $s>0$ and $s<0$. (With the common  point $X=s=0$;
note that at $s=0$ the left hand side of (14) factors as  $(x^2-1)^2(x^2+k^2)$.)  In variables $(X,t)$ it may be described 
somewhat more explicitly as
$$t\ge t_{min}(X)\ge 0$$
where  $t_{min}(X)=\frac{1}{2}X^2+O(X^4)$, but  the function $t_{min}$ depends on the choice between $s=\sqrt{t}$ and $s=-\sqrt{t}$.

The rest of the construction  is similar to \cite{Gav}. We introduce the function $y$ by 
$$dy=\pm\frac{dt-Fdx}{\sqrt{G}}; \eqno{(15)}$$
the form on the right hand side  of (15) is closed because of (8) and exact because the domain  may be chosen  simply connected. 
We may assume that $y=0$ for $t=t_{min}$, which allows us to extend $t$ to a function $t(x,y)=t(x,-y)$ analytic near $(x,y)=(1,0)$
except for the point itself. (Not all of this is immediately  obvious, but the argument  is the same as in \cite[Lemma 3]{Gav}.) 

\subsection{The modification}

As explained in \cite[\S 3]{Gav}, a localizable (i.e. satisfying  $u\cdot\nabla p=0$)  Euler flow $(u,p)$  can be modified 
to obtain another Euler flow,   
$$\widetilde{u}=\omega u,\,d\widetilde{p}=\omega^2 \,dp.$$
Choosing a smooth function  $\omega$ such that  $\omega(t)=0$ for $t\not\in [\varepsilon, 2\varepsilon]$, 
we can obtain  a smooth Euler flow $\widetilde{u}$ with support near the  helix. (This is why we call such a flow  ``localizable'':  
it allows  modifications which can  reduce  its  support.)

\subsection{A comparison with the circle}

It is possible to consider the Euler flow $u$ constructed  in \cite{Gav} as a degenerate  case  corresponding to $k=0$ (when our  helix turns  into a circle). In terms of \cite{Gav} for our flow we have  $R=1$ and 
$$\alpha=4t,\,H(\alpha)=16h^2(t).$$
If  $k=0$ then  
$$\frac{\dd}{\dd x}t=F=\frac{1}{2}x(x^2-c);$$
comparing this with 
$$\frac{\dd}{\dd x}\alpha=2x^3-2x\psi(\alpha)$$
in \cite{Gav} we must conclude that $c(t)=\psi(\alpha)$. Now  (12) become   
$$\frac{dh^2}{dt}=c+\frac{36t^2}{h^2-3tc},\,\,\,\frac{dc}{dt}=\frac{6t}{h^2-3tc};$$
excluding $h^2$ from this system we have a second order equation 
$$6tc''+3t(c')^3-2c(c')^2-6c'=0,$$
the same as in \cite[Lemma 1]{Gav}. 

It should be noted that for $k=0$ the function $t=t(x,y)=\frac{1}{4}\alpha$ is analytic at the point $(1,0)$ while for $k>0$ 
it is  not, although this difference  cannot be seen from just the main term  of  the asymptotic, 
$$t(x,y)\sim \frac{(x-1)^2}{2}+ \frac{y^2}{2(1+k^2)}.$$
This fact is related to the choice between two solutions mentioned above, which  actually correspond   to different analytic branches of this function. 
(The  branching  curve of $t$  in the complex $(x,y)$  plane  has only one real point, so outside of it this  function is real analytic.)

\section{Some observations}

\subsection{Beltrami flows}

The author would like to point out that  an Euler flow satisfying the condition\footnote{Or  any relation of the form $|u|^2=f(p)$ for that matter.}
$|u|^2=3p$ may be interpreted as a special case of a Beltrami flow. Indeed,  a modification 
$$\widetilde{u}=\omega u,\,d\widetilde{p}=\omega^2 \,dp$$
with $\omega=p^{-\frac{5}{6}}$ satisfies  $\frac{1}{2}|\widetilde{u}|^2+\widetilde{p}=0,$
so by the  known identity \cite[\S 2]{LL}
$$\frac{1}{2}\nabla |u|^2=u\times\curl u+(u\cdot\nabla) u \eqno{(16)}$$
we have $\widetilde{u}\times \curl \widetilde{u}=0$ i.e. $\curl \widetilde{u}=\lambda  \widetilde{u}$
for some function  $\lambda$. This Beltrami flow is localizable  
because   $|u|^2=3p$  implies $u\cdot \nabla p=0$. Conversely, given a localizable  Beltrami flow $ \widetilde{u}$
one can modify it to obtain a solution with $|u|^2=3p$. From the theoretical perspective  Beltrami flows
are more convenient  to consider, so we will take this point of view in this section.

\subsection{The Grad-Shafranov equation}

Let 
$$\xi=\dd_\varphi+k\dd_z=\rho e_\varphi+ke_z$$
be the  Killing  vector field from above. Under close examination, the field 
$$a=\frac{\xi}{|\xi|^2}$$
is another  Beltrami flow, 
$$\ddiv a=0,\,\curl a=2k|\xi|^{-2}a.$$ 
The Beltrami modification $\widetilde{u}$ of the field $u$ we have constructed in \S 4 may be written in the form
$$\widetilde{u}=a\times\nabla \psi+\chi a \eqno{(18)}$$
where  $\psi=t^{\frac{1}{6}},\,\chi=\frac{1}{6}t^{-\frac{5}{6}}h(t)$. 

 Now we may change the perspective and ask   the following question: if  $\widetilde{u}$ is some  field given by (18) with  $\xi(\psi)=\xi(\chi)=0$,
what additional  conditions this  two functions must  satisfy to make it  a Beltrami flow? 
Note that  we have $\ddiv \widetilde{u}=0$ automatically, so this  is a  question  about the vorticity. 
Using the formula\footnote{Where $[\cdot,\cdot]$ is the Lie bracket.}  \cite[Ch II, \S 1]{AK}
$$\curl(A\times B)=(\ddiv B)A-(\ddiv A)B-[A,B], \eqno{(19)}$$
we have 
$$\curl \widetilde{u}=(\Delta\psi +2k\chi|\xi|^{-2})a-[a,\nabla\psi]+\nabla\chi\times a.$$
Note that $\xi$ is a generator of  isometry and $\xi(\psi)=0$. It follows that $[\xi,\nabla\psi]=0$ hence 
$[a,\nabla\psi]=2(\nabla \log |\xi|,\nabla\psi)a$ and 
$$\curl \widetilde{u}=(\Delta\psi +2k\chi|\xi|^{-2}-2(\nabla \log |\xi|,\nabla\psi))a+\nabla\chi\times a.$$
Assuming that 
$$\curl \widetilde{u}=\lambda\widetilde{u}=\lambda\chi a-\lambda \nabla \psi\times a$$
and taking into account that both gradients are orthogonal to $a$, we must conclude that 
$$d\chi=-\lambda d\psi,$$
which essentially  implies $\chi=h(\psi),\,\lambda=-h'(\psi)$ (for some function $h$).
Then the  condition  can be written in the form
$$\Delta\psi -2(\nabla \log |\xi|,\nabla\psi)+2k |\xi|^{-2} h+hh'=0.\eqno({20})$$
In the case  $k=0$ we have $|\xi|=\rho$ and the equation  becomes 
$$\Delta\psi -\frac{2}{\rho}\frac{\dd\psi}{\dd\rho}+hh'=0, \eqno{(21)}$$
which  is known as   (a special case of)  the  Grad-Shafranov equation. 

One can see  that our  construction  was actually  built on a Killing  field.  If we drop the assumption that $\xi$ 
is a generator of isometry  then, apparently, we have no means to control the Lie bracket and the whole construction falls apart. 
It  looks like  the axial or helical  symmetry of a  flow  was  not merely a simplification to make the calculation
easy  but  is necessary to make the ends meet. If there are any localizable Euler  flows 
which are not symmetric  the author does not really  know. 

\subsection{Special Beltrami flows on Riemannian manifolds?}

Unfortunately,  the Euclidean  space  has  no one-parameter isometry  groups  besides what  we have already considered. However, we may take a more broad view and ask about possible generalizations of the above construction to Riemannian manifolds. The Euler equation on a 
Riemannian manifold  is  the same as in the Euclidean space, except  $\nabla u$ must now  be interpreted as the  covariant derivative  of $u$ with respect to the Levi-Civita  connection \cite{AK}.  A Belirami flow  on  an oriented  manifold of dimension three is defined as usual, it is a vector field  $u$ satisfying
$$\curl u=\lambda u,\,\,\,\ddiv u=0.$$ 
The fact that such a flow obeys the Euler equation\footnote{With $p=-\frac{1}{2}|u|^2+\text{const}$.}  follows from the  formula (16).  
We assume forth  that our Riemannian manifold has a Killing vector field $\xi$. In this case (\S 5.4)
$$\xi\times\curl\xi=\nabla |\xi|^2,\eqno{(22)}$$ 
and   the field $a=|\xi|^{-2}\xi$ is again a Beltrami flow. Indeed, $\ddiv a=0$ because $\xi(|\xi|^{-2})=0$ and $\ddiv\xi=0$, and 
$$a\times\curl a=|\xi|^{-4}\xi\times\curl\xi+|\xi|^{-2} \xi\times (\nabla |\xi|^{-2}\times\xi)=|\xi|^{-4}\nabla |\xi|^2-|\xi|^{-6}\xi\times (\nabla |\xi|^{2}\times\xi) =0,$$
because the vectors $\xi$ and $\nabla |\xi|^2$ are orthogonal. Thus, we have 
$$\ddiv a=0,\,\curl a=\mu a$$
for some function $\mu$. 

We may again  try to construct a  Belirami flow of the form 
$$u=a\times\nabla \psi+\chi a,$$
assuming that  $\xi(\psi)=\xi(\chi)=0$. In this case $\ddiv u=0$ and
$$\curl u=\Delta\psi a-[a,\nabla \psi]+\chi\mu a+\nabla\chi\times a.$$   
Repeating the computation from \S  5.3,  we have  $\chi=h(\psi)$ and
$$\Delta\psi -2(\nabla \log |\xi|,\nabla\psi)+h(h'+\mu)=0,\eqno({23})$$
which may be considered a generalization of the Grad-Shafranov equation (21).

The way to make this Beltrami flow localizable  is to assume that $|u|^2$ depends on  $\psi$.
We obviously have  $|u|^2=|\xi|^{-2}(|\nabla \psi|^2+h^2(\psi)),$ so this condition means 
$$|\nabla \psi|^2=|\xi|^{2} f(\psi)-h^2(\psi)\eqno{(24)}$$
for some function $f$. All the variables in (23,24) are invariant under isometries  generated by  $\xi$, so this 
is, in fact,   a PDE  in the (two-dimensional) space of orbits rather then in the original  manifold. 
The problem is that  it is overdetermined and does not seem easy to handle.
To overcome this obstacle in our  special case  we have \emph{de facto}  introduced a somewhat contrived vector field $U(\psi)$ depending on $\psi$ as a parameter, and then showed that both (23) and (24) follow from the same  equation
$$\nabla\psi=U(\psi).\eqno{(25)}$$
However, it may be difficult (if possible at all)  to pull off  the same trick in the general case. 

\subsection{Proofs of some formulas} 

We will prove here  vector calculus  formulas  (16, 19, 22)  used in the last section. All of this  proofs are very simple, but the formulas are important 
and the author has yet to  see  them derived properly in the literature. So, he decided to  write them down for the sake of a
reader's convenience. 

We are dealing with  an oriented Riemannian manifold of dimension three. Let $(x_1, x_2,x_3)$ 
be  local coordinates  and $\dd_i=\frac{\dd}{\dd x^i}$ be the corresponding   vector fields;
naturally, we assume that the  frame $(\dd_1, \dd_2, \dd_3)$ agrees with the orientation.  For a vector field $u= u^i\dd_i$ 
we have the following basic  formulas 
$$(\curl u)^i=\epsilon^{ijk} \nabla_j u_{k},\, [(u\cdot\nabla) u]^i=u^j \nabla_j u^i,\,\,\,(u\times v)_i=\epsilon_{ijk}u^jv^k,$$
where $\epsilon$ is the Levi-Civita tensor\footnote{Not the  Levi-Civita \emph{symbol}. In this notation,  $\epsilon_{123}=\sqrt{|g|}$ and $\epsilon^{123}=1/\sqrt{|g|}$.}.

Unfortunately, from the presentation in  the book  \cite[Ch. II, \S 1]{AK}  it is unclear  if  the  formula (16)  was only meant   for the Euclidean space  or is valid regardless of  the metric. However, it is not difficult to see that  the latter is the case (there are no higher order derivatives 
of the metric, so  curvature terms do not appear). Indeed, we have   
$$[u\times \curl u]^i+[(u\cdot\nabla) u]^i=\epsilon_{klj}\epsilon^{k\alpha\beta}g^{il} u^j\nabla_\alpha u_{\beta}+u^j \nabla_j u^i=$$
$$=u^j(\nabla^i u_j - \nabla_j u^i )+u^j \nabla_j u^i=\frac{1}{2}(\nabla |u|^2)^i.$$

The  formula  for the curl  of a cross  product  is  invariably omitted from  vector calculus textbooks because of the  Lie bracket 
(which, apparently,  is considered  inappropriate for undergraduates).   If   $C=A\times B$, then 
$$\nabla_j C_{i}=\epsilon_{kli}(B^l\nabla_j A^k+A^k\nabla_j B^l)$$
because $\nabla \epsilon=0$. Thus
$$(\curl C)^i=\epsilon^{ijk}\nabla_j C_{k}=\epsilon^{ijk}  \epsilon_{mlk}(B^l\nabla_j A^m+A^m\nabla_j B^l)=$$
$$(B^j\nabla_j A^i+A^i\nabla_j B^j)-(B^i\nabla_j A^j+A^j\nabla_j B^i)=A^i\ddiv  B-B^i\ddiv A-[A,B]^i.$$

A proof of the last one is not any more complicated. If  $\xi$ is a Killing vector, then 
$\nabla_j\xi_{i}+\nabla_i\xi_{j}=0$ by definition, hence  
$$\frac{1}{2}(\nabla |\xi|^2)^i=\xi^j \nabla^i\xi_{j}=-\xi^j \nabla_j\xi^i=-[(\xi\cdot\nabla)\xi]^i.$$
Comparing this equality with (16), we have
$$\xi\times\curl\xi=\nabla |\xi|^2.$$

{}

\end{document}